\documentclass{amsproc}
\usepackage{amssymb}
\usepackage{graphicx}
\usepackage{amscd}
\usepackage{amsmath}
\theoremstyle{plain}

\newtheorem{lemma}{Lemma}

\newtheorem{proposition}{Proposition}
\newtheorem{remark}{Remark}

\newtheorem{theorem}{Theorem}
\newtheorem{assumption}{Assumption}
\numberwithin{equation}{section}

\newcommand{\N}{\mathbb{N}}
\newcommand{\p}{\partial}

\begin{document}

\title{On some applications of the Boundary Control method to spectral estimation and inverse problems }
\author{S. A. Avdonin}
\address{Department of Mathematics and Statistics, University of Alaska,
Fairbanks, AK 99775-6660, USA.} \email{s.avdonin@alaska.edu}
\author{ A. S. Mikhaylov}
\address{St. Petersburg   Department   of   V.A. Steklov    Institute   of   Mathematics
of   the   Russian   Academy   of   Sciences, 7, Fontanka, 191023
St. Petersburg, Russia and St.-Petersburg State University, Faculty
of Mathematics and Mechanics, Universitetskii av. 28, Petrodvorec,
198504 St-Petersburg, Russia.} \email{mikhaylov@pdmi.ras.ru}

\author{ V. S. Mikhaylov}
\address{St.Petersburg   Department   of   V.A.Steklov    Institute   of   Mathematics
of   the   Russian   Academy   of   Sciences, 7, Fontanka, 191023
St. Petersburg, Russia and St.-Petersburg State University, Faculty
of Physics, Ulyanovskaya str. 3, Petrodvorec, 198504
St-Petersburg, Russia.} \email{vsmikhaylov@pdmi.ras.ru}
\date{November 26, 2014}

\maketitle

\noindent {\bf Abstract.} We consider applications of the Boundary
Control (BC) method to generalized spectral estimation problems
and to inverse source problems. We derive the equations of the BC
method for this problems and show that solvability of this
equations crucially depends on the controllability properties of
the corresponding dynamical system and properties of corresponding
families of exponentials.

\section{Introduction.}

The classical spectral estimation problem consists of the recovery
of the coefficients $a_n$, $\lambda_k$, $k = 1,\ldots,N, \; N \in
\mathbb{N},$ of a signal
\begin{equation*}
s(t)=\sum_{n=1}^N a_ke^{\lambda_k t},\quad t\geqslant 0
\end{equation*}
from the given observations $s(j)$, $j=0,\ldots, 2N-1,$ where the
coefficients $a_k,$ $\lambda_k$ may be arbitrary complex numbers.
The literature describing variuos methods for solving the spectral
estimation problem is very extensive: see for example the list of
references in \cite{AB,ABN}. In these papers a new approach to
this problem was proposed: a signal $s(t)$ was treated as a kernel
of certain convolution operator corresponding to an input-output
map for some linear discrete-time dynamical system. While the
system realized from the input-output map is not unique, the
coefficients $a_n$ and $\lambda_n$ can be determined uniquely
using the non-selfadjoint version of the boundary control method
\cite{ABe}.

In \cite{AGM,AM2} this approach has been generalized to
infinite-dimensional case: more precisely, the problem of the
recovering the coefficients $a_k, \lambda_k\in \mathbb{C}$, $k\in
\mathbf{N},$ of the given signal
\begin{equation}
\label{signal} S(t)=\sum_{k=1}^\infty a_k(t)e^{\lambda_kt},\quad
t\in (0,2T),
\end{equation}
from the given data $S\in L_2(0,2T)$ was considered. In \cite{AGM}
the case $a_k\in \mathbb{C}$ has been treated, in \cite{AM2} the
case when for each $k,$ $a_k(t)=\sum_{i=0}^{L_k-1}a^i_kt^i$ are
polynomials of the order $L_k-1$ with complex valued coefficients
$a^i_k$ was studied.

Recently it was observed \cite{AMR,MM} that the results of
\cite{AGM,AM2} are closely related to the dynamical inverse source
problem: let $H$ be a Hilbert space, $A$ be an operator in $H$
with the domain $D(A)$, Y be another Hilbert space, $O:H\supset
D(O)\mapsto Y$ be an observation operator (see \cite{TW}). Given
the dynamical system in $H$:
\begin{equation}
\label{E_D}
\left\{
\begin{array}l
u_t-Au=0,\quad t>0,\\
u(0)=a,
\end{array}
\right.
\end{equation}
we denote by $u^a$ its solution, and by $y(t):=(Ou^a)(t)$ the
observation (output of this system). The operator that realize the
correspondence $a\mapsto (Ou^a)(t)$ is called \emph{observation}
operator $\mathbb{O}^T: H\mapsto L_2(0,T;Y)$. We fix some $T>0$
and assume that $y(t)\in L_2(0,T;Y)$. One can pose the following
questions: what information on the operator $A$ could be recovered
from the observation $y(t)$? We mention works on the
multidimensional inverse problems for the Schr\"odinger, heat and
wave equations by one measurement, concerning this subject. Some
of the results (for the Schr\"odinger equation) are given in
\cite{BP,MOR,AMR}. To answer this question in the abstract
setting, in \cite{MM} the authors derived the version of the
BC-method equations under the condition that $A$ is self-adjoint
and $Y=\mathbb{R}$. In the present paper we address the same
question without the assumption about selfajointness of $A$ . The
possibility of recovering the spectral data from the dynamical one
is well-known for the dynamical system with boundary control
\cite{B01,B03}. We extend this ideas to the case of the dual
(observation) system.

The solvability of the BC-method equations for the spectral
estimation problem critically depends on the properties of
corresponding exponential family. The solvability of the BC-method
equations for the system (\ref{E_D}) depends on the
controllability properties of the dual system. We point out the
close relation between these two problems: they both leads to
essentially the same equations (see section 4 for applications),
and conditions of the solvability of these equations are the seme
(on the connections between the controllability of a dynamical
systems and properties of exponential families see \cite{AI}).

In the second section we outline the solution of the spectral
estimation problem in infinite dimensional spaces (see \cite{AM2}
for details). In the third section we derive the equations of the
BC-method for the problem (\ref{E_D}) extending the results of
\cite{MM} to the case of non self-adjoint operator. Also we answer
the question on the extension of the observation
$y(t)=(\mathbb{O}a)(t)$. The last section is devoted to the
applications to inverse problem by one measurement for the
Schr\"odinger equation on the interval and to the problem of the
extension of the inverse data for the first order hyperbolic
system on the interval, see also \cite{AGM,AM1,AM2,AMR}.

\section{The spectral estimation problem in infinite dimensional spaces.}

The problem is set up in the following way: given the signal
(\ref{signal}), $S\in L_2(0,2T)$, for $T>0$, to recover the
coefficients $a_k(t)$, $\lambda_k$, $k\in \mathbb{N}.$ Below we
outline the procedure of recovering unknown parameters, for the
details see \cite{AM2}.

We consider the dynamical systems in a complex Hilbert space $H$:
\begin{eqnarray}
\label{dyn_sys} \dot x(t)=Ax(t)+bf(t), \quad t\in (0,T),\quad
x(0)=0.\\
\dot y(t)=A^*y(t)+dg(t), \quad t\in (0,T),\quad
y(0)=0,\label{dyn_sys_adj}
\end{eqnarray}
Here $b,d\in H$, $f,g\in L_2(0,T),$ and we assume that the
spectrum of the operator $A$, $\{\lambda_k\}_{k=1}^\infty$ is not
simple. We denote the algebraic multiplicity of $\lambda_k$  by
$L_k,$ $k \in \mathbb{N},$ and assume also that the set of all
root vectors $\{\phi_k^i\},$ $i=1,\ldots, L_k,$  $k\in
\mathbb{N},$ forms a Riesz basis in $H$. Here the vectors from the
chain $\{\phi^i_k\}_{i=1}^{L_k}$, $k \in \mathbb{N},$ satisfy the
equations
\begin{equation*}
\left(A-\lambda_k\right)\phi_k^1=0,\quad
\left(A-\lambda_k\right)\phi_k^i=\phi_k^{i-1},\,\, 2\leqslant
i\leqslant L_k.
\end{equation*}
The spectrum of $A^*$ is $\{\overline\lambda_k\}_{k=1}^\infty$ and
the root vectors $\{\psi_{k}^{i}\}$, $i=1,\ldots, L_k,$  $k\in
\mathbb{N},$  also form a Riesz basis in $H$ and satisfy the
equations
\begin{equation*}
\left(A^*-\overline\lambda_k\right)\psi_{k}^{L_k}=0,\quad
\left(A^*-\overline\lambda_k\right)\psi_{k}^i=\psi_{k}^{i+1},\,\,
1\leqslant i\leqslant L_k-1.
\end{equation*}
Moreover, the root vectors of $A$ and $A^*$ are normalized in
accordance with
\begin{eqnarray*}
\left<\phi_k^i,\psi_{l}^j\right>=0, \text{ if $k\not=l$ or
$i\not= j$};\\
\left<\phi_k^i,\psi_{k}^i\right>=1,\,\,\,i=1,\ldots,L_k, \, k\in
\mathbb{N}.
\end{eqnarray*}
We  consider $f$ and $g$ as the inputs of  the systems
(\ref{dyn_sys}) and (\ref{dyn_sys_adj}) and define the outputs $z$
and $w$ by the formulas
\begin{equation*}
z(t)=\left<x(t),d\right> ,\quad w(t)=\left<y(t),b\right>.
\end{equation*}
Suppose that $b=\sum_{k=1}^\infty\sum_{i=1}^{L_k}b_k^i\phi_{k}^i$,
$d=\sum_{k=1}^\infty\sum_{i=1}^{L_k}d_k^i\psi_{k}^i.$ Looking for
the solution to (\ref{dyn_sys}) in the form
$x(t)=\sum_{k=1}^\infty\sum_{i=1}^{L_k}c_k^i(t)\phi_{k}^i$, we
arrive at the following representation for the output
\begin{equation*}
z(t)=\left<x(t),d\right>=\sum_{k=1}^\infty\sum_{i=1}^{L_k}c_k^i(t)d_k^i=\int_0^tr(t-\tau)f(\tau)\,d\tau,
\end{equation*}
where the {\sl response function} $r(t)$ is defined as
\begin{equation}
\label{r_representation} r(t)=\sum_{k=1}^\infty e^{\lambda_k
t}\left[a_k^1+a_k^2t+a_k^3\frac{t^2}{2}+\ldots
+a_k^{L_k-1}\frac{t^{L_k-2}}{(L_k-2)!}+a_k^{L_k}\frac{t^{L_k-1}}{(L_k-1)!}\right],
\end{equation}
with $a_k^j$ being defined
\begin{equation}
\label{coeff_repr} a_k^j=\sum_{i=j}^{L_k}b_k^i d_k^{i-j+1},\quad
 j=1,\ldots, L_k, \  k\in
\mathbb{N}.
\end{equation}
It is important to note  that
$r(t)$ has the form of the series in (\ref{signal}).

Analogously, looking for the solution of (\ref{dyn_sys_adj}) in
the form
\begin{equation*}
y(t)=\sum_{k=1}^\infty\sum_{i=1}^{L_k}h_k^i(t)\psi_{k}^i,
\end{equation*}
we arrive at
\begin{equation*}
w(t)=\left<y(t),b\right>=\sum_{k=1}^\infty\sum_{i=1}^{L_k}h_k^i(t)b_k^i=
\int_0^t \overline{r(t-\tau)}g(\tau)\,d\tau.
\end{equation*}
We introduce the {\sl connecting operator} $C^T: L_2(0,T)\mapsto
L_2(0,T)$ defined through its bilinear form by the formula:
\begin{equation*}
\left<C^Tf,g\right>=\left<x(T),y(T)\right>.
\end{equation*}
In \cite{AM2} the representation for $C^T$ was obtained:
\begin{lemma}
The connecting operator $C^T$ has a representation
\begin{equation*}
(C^Tf)(t)=\int_0^T r(2T-t-\tau)f(\tau)\,d\tau.
\end{equation*}
\end{lemma}
We assume that the systems (\ref{dyn_sys}), (\ref{dyn_sys_adj})
are spectrally controllable in time $T$. This means that, for any
$i \in \{1,\ldots,L_k\}$ and  any $k \in \mathbb{N},$  there exist
$f_k^i,\, g_k^i \in H^1_0(0,T)$,
 such that $x^{f_k^i}(T)=\phi_k^i$, $y^{g_k^i}(T)=\psi_{k}^i$.
Using ideas of the BC method \cite{B07}, we are able to extract
the spectral data, $\left\{\lambda_k, a_k^j\right\},$ $j=1,\ldots,
L_k, \  k\in \mathbb{N}$, from the dynamical one, $r(t), \; t \in
(0,2T),$ (see \cite{AGM,AM2} for more details):
\begin{proposition}
The set $\lambda_k,$ $f_k^i$,  $i=1,\ldots,L_k,$ $k \in \N ,$ are
eigenvalues and root vectors of the following generalized
eigenvalue problem in $L_2(0,T)$:
\begin{equation}
\label{gen_ev_prob} \int_0^T\left(r'(2T-t-\tau)-\lambda
r(2T-t-\tau)\right)f(\tau)\,d\tau=0.
\end{equation}
The set $\overline\lambda_k,$ $g_k^i$, $k=1,\ldots\infty$,
$i=1,\ldots,L_k$ are eigenvalues and root vectors of the
generalized eigenvalue problem in $L_2(0,T)$:
\begin{equation}
\label{gen_ev_prob_adj} \int_0^T\left(\overline{r'(2T-t-\tau)}
-\lambda \overline{r(2T-t-\tau)}\right)g(\tau)\,d\tau=0.
\end{equation}
\end{proposition}

Now we describe the algorithm of recovering $a_k^1,\ldots
a_k^{L_k}$, $k \in \N$ (see the representation
(\ref{r_representation})). We normalize the solutions to
(\ref{gen_ev_prob}), (\ref{gen_ev_prob_adj}) by the rule
\begin{equation}
\label{norming} \left<C^T\widetilde f_k^i,\widetilde
g_k^i\right>=1.
\end{equation}
and define
\begin{eqnarray}
\widetilde b_k^i=\left<y^{\widetilde
g_k^i}(T),b\right>=\int_0^T\overline r(T-\tau)\widetilde
g_k^i(\tau)\,d\tau,\label{b_i}\\
\widetilde d_k^i=\left<x^{\widetilde f_k^i}(T),d\right>=\int_0^T
r(T-\tau)\widetilde f_k^i(\tau)\,d\tau.\label{d_i}
\end{eqnarray}
Then (see (\ref{coeff_repr}))
\begin{equation}
\label{a_k_1} a_k^1=\sum_{i=1}^{L_k}\widetilde b_k^i\widetilde
d_k^i.
\end{equation}

Denote by $\partial$ and $I$ the operator of differentiation and
the identity operator in $L_2(0,T).$ We normalize the solutions to
(\ref{gen_ev_prob}), (\ref{gen_ev_prob_adj}) (for $i>l$) by the
rule
\begin{equation}
\label{norming_l}
\left<\left[C^T\left(\partial-\lambda_kI\right)\right]^l\widehat
f_k^i,\widehat g_k^{i-l}\right>=1,
\end{equation}
we define $\widehat b_k^i,$ $\widehat d_k^i$ by (\ref{b_i}),
(\ref{d_i}) and evaluate
\begin{equation}
\label{a_k_l} a_k^l=\sum_{i=l}^{L_k}\widehat b_k^i \widehat
d_k^{i-l+1},\quad l=2,\ldots,L_k.
\end{equation}
We conclude this section with the algorithm for solving the
spectral estimation problem: suppose that we are given with the
function $r\in L_2(0,2T)$ of the form (\ref{r_representation}) and
the family
$\bigcup_{k=1}^\infty\{e^{\lambda_kt},\ldots,t^{L_k-1}e^{\lambda_kt}\}$
is minimal in $L_2(0,T)$. Then to recover $\lambda_k$, $L_k$ and
coefficients of polynomials, one should follow the

{\bf Algorithm}
\begin{itemize}
\item[a)] solve generalized eigenvalue problems
(\ref{gen_ev_prob}), (\ref{gen_ev_prob_adj}) to find $\lambda_k$,
$L_k$ and non-normalized controls.

\item[b)] Normalize $\widetilde f_k^i,$ $\widetilde g_k^i$ by
(\ref{norming}), define $\widetilde b_k^i,$ $\widetilde d_k^i$ by
(\ref{b_i}), (\ref{d_i}) to recover $a_k^1$ by (\ref{a_k_1}).

\item[c)] Normalize $\widehat f_k^i,$ $\widehat g_k^{i-l}$ by
(\ref{norming_l}), define $\widehat b_k^i,$ $\widehat d_k^i$ by
(\ref{b_i}), (\ref{d_i}) to recover $a_k^l$ by (\ref{a_k_l}),
$l=2,\ldots,L_k-1$.

\end{itemize}

\section{Equations of the BC method.}

Let us denote by $A^*$ the operator adjoint to $A$ and $B:=O^*$,
$B: Y\mapsto H$. Along with the system (\ref{E_D}) we consider the
following dynamical control system:
\begin{equation}
\label{E_AD}
\left\{
\begin{array}l
v_t+A^*v=Bf,\quad t<T,\\
v(T)=0,
\end{array}
\right.
\end{equation}
and denote its solution by $v^f$. The reason we consider the
system (\ref{E_AD}) backward in time is that it is adjoint to
(\ref{E_D}) (see \cite{AI, MM}).

For every $0\leqslant s<T$ we
introduce the \emph{control} operator by $W^sf:=v^f(s).$ It is easy
 to check that $-W^0$ is adjoint to $\mathbb{O}^T.$  Indeed,
taking $f\in L_2(0,T;Y),$ $a\in H$ we show \cite{MM} that
\begin{equation}
\label{Adj}
\int_0^T\left(f,\mathbb{O}a\right)_Y=-\left(W^0f,a\right)_H,
\end{equation}
here $\mathbb{O}a=\left(Ou^a\right)(t).$ Due to the arbitrariness
of $f$ and $a$, the last equality is equivalent to
$\left(\mathbb{O}^T\right)^*=-W^0$.

\noindent We assume that the operator $A$ satisfies the following
assumptions:
\begin{assumption}
\label{A1}
\begin{itemize}
\item[a)] The spectrum of the operator $A$,
$\{\lambda_k\}_{k=1}^\infty$ consists of the eigenvalues
$\lambda_k$  with algebraic multiplicity $L_k,$ $k \in
\mathbb{N},$ and the set of all root vectors $\{\phi_k^i\},$
$i=1,\ldots, L_k,$  $k\in \mathbb{N},$ form a Riesz basis in $H$.
Here the vectors from the chain $\{\phi^i_k\}_{i=1}^{L_k}$, $k \in
\mathbb{N},$ satisfy the equations
\begin{equation*}
\left(A-\lambda_k\right)\phi_k^1=0,\quad
\left(A-\lambda_k\right)\phi_k^i=\phi_k^{i-1},\,\, 2\leqslant
i\leqslant L_k.
\end{equation*}
The root vectors of $A^*$, $\{\psi_k^i\},$ $i=1,\ldots, L_k,$
$k\in \mathbb{N},$ form a Riesz basis in $H$ and satisfy:
\begin{equation*}
\left(A^*-\overline{\lambda}_k\right)\psi_k^{L_k}=0,\quad
\left(A^*-\overline{\lambda}_k\right)\psi_k^i=\psi_k^{i+1},\,\,
1\leqslant i\leqslant L_k-1.
\end{equation*}
\item[b)] The system (\ref{E_AD}) is spectrally controllable in
time $T$: i.e. there exists the controls $f_k^i\in H^1_0(0,T;Y)$
such that $W^0f^i_k=\psi^i_k$, for  $i=1,\ldots, L_k$, $k \in
\mathbb{N}$.
\end{itemize}
\end{assumption}

We say that the vector $a$ is \emph{generic} if  its Fourier
representation in the basis $\{\phi_k^i\}_{k=1}^\infty$,
$a=\sum_{k=1}^\infty\sum_{i=1}^{L_k} a_k^i\phi_k^i$, is such that
$a_k^i\not=0$ for all $k,i$. We assume that the controls from the
Assumption \ref{A1} are extended by zero outside the interval
$(0,T)$. Now we are ready to formulate

\begin{theorem}
\label{teor}
If $A$ satisfies Assumption \ref{A1},  $Y=\mathbb{R},$ and the source $a$ is
generic, then the spectrum of $A$ and controls $f_k^i$ are the
spectrum and the root vectors of the following generalized
spectral problem:
\begin{equation}
\label{M_eqn}
\int_0^{2T}\left(\dot{(Oa)}(t)-\lambda_k(Oa)(t),f_k(t-T+\tau)\right)_Y\,dt=0,\quad
0<\tau<T.
\end{equation}
\end{theorem}
Here by dot we
denote the differentiation with respect to $t$.
\begin{proof}
 We
denote by $\{\widetilde f_k^i\}$ the set of controls which satisfy
$W^0\widetilde f_k^i=\psi_k^i$. By $\{f_k^i\}$ we denote the set
of shifted controls: $f_k^i(t)=\widetilde f_k^i(t-T)$. Thus the
controls $f_k^i$ acts on the time interval $(T,2T)$.  Let us fix some
 $i\in 1,\ldots,L_k$, $k\in \mathbb{N},$ $\tau\in (0,T)$ and
consider $W^0\left(\dot f_k^i(\cdot+\tau)\right)$:
\begin{equation}
\label{w_1} W^0\left(\dot f_k^i(\cdot+\tau)\right)=v^{\dot
f_k^i(\cdot+\tau)}(0)=v_t^{
f_k^i(\cdot+\tau)}(0)=\left(Bf_k^i(\cdot+\tau)\right)(0)-A^*v^{f_k^i(\cdot+\tau)}(0).
\end{equation}
Since $f_k^i\in H^1_0(T,2T,Y)$,
$\left(Bf_k^i(\cdot+\tau)\right)(0)=0.$ The second term in the
right hand side of $(\ref{w_1})$ could be evaluated using the
following reasons. The function $v^{f_k^i}$ solves:
\begin{eqnarray*}
&v^{f_k^i(\cdot+\tau)}_t+A^*v^{f_k^i(\cdot+\tau)}=0,\quad 0\leqslant t\leqslant T-\tau,\\
&v^{f_k^i(\cdot+\tau)}(T-\tau)=\psi_k^i.
\end{eqnarray*}
We are looking for the solution in the form
$v^{f_k^i(\cdot+\tau)}(t)=\sum_{j=1}^{L_k} c_k^j(t)\psi_k^j$ then
$c_k^j$ satisfy boundary condition $c_k^j(0)=\delta_{ij}$ and
equation:
\begin{align*}
&\frac{d}{dt}c_k^1+\overline{\lambda}_kc_k^1=0,\\
&\frac{d}{dt}c_k^j+\overline{\lambda}_kc_k^j+c_k^{j-1}=0,\quad
j=2,\ldots, L_k.
\end{align*}
Solving this system we obtain the following expansion
\begin{equation}
\label{v_f}
v^{f_k^i(\cdot+\tau)}(t)=\sum_{j=i}^{L_k}\frac{(T-\tau-t)^{j-i}}{(j-i)!}e^{\overline{\lambda}_k(T-\tau-t)}\psi^j_k
\end{equation}
Evaluating $A^*v^{f_k^i(\cdot+\tau)}(0),$ making use of (\ref{v_f})
and properties of the root vectors, we arrive at:
\begin{align*}
&A^*v^{f_k^{L_k}(\cdot+\tau)}(0)=\overline{\lambda}_kv^{f_k^{L_k}(\cdot+\tau)}(0),\\
&A^*v^{f_k^i(\cdot+\tau)}(0)=\overline{\lambda}_kv^{f_k^i(\cdot+\tau)}(0)+v^{f_k^{i+1}(\cdot+\tau)}(0),\,\,
i<L_k.
\end{align*}
Then continuing (\ref{w_1}), we obtain:
\begin{eqnarray}
 W^0\left(\dot
 f_k^{L_k}(\cdot+\tau)\right)=-A^*v^{f_k^{L_k}(\cdot+\tau)}(0)=-\overline{\lambda}_kW^0f_k^{L_k},\label{w_31}\\
W^0\left(\dot
 f_k^{i}(\cdot+\tau)\right)=-\overline{\lambda}_kW^0f_k^{i}-\overline{\lambda}_kW^0f_k^{i+1},\,\,i<L_k.\label{w_32}
\end{eqnarray}

Integrating by parts and taking into account that
$f_k^i(0)=f_k^i(T)=0$ for $i=1,\ldots,L_k,$ we get:
\begin{eqnarray}
\int_0^{2T}\left((Oa)(t),\dot f_k^i(t+\tau)\right)_Y\,dt=-\int_{0}^{2T}\left(\dot{(Oa)}(t),f_k^i(t+\tau)\right)_Y\,dt\notag\\
+\left(\dot{(Oa)}(t+\tau),f_k^i(t)\right)_Y\Bigr|_{t=0}^{t=2T}=-\int_{0}^{2T}\left(\dot{(Oa)}(t),f_k^i(t+\tau)\right)_Y\,dt\label{OD}
\end{eqnarray}
One the other hand, using the duality between $W^0$ and
$\mathbb{O}^T$ and (\ref{w_31}), (\ref{w_32}), we have for
$i=L_k$:
\begin{eqnarray}
\int_0^{2T}\left((Oa)(t),\dot
f_k^{L_k}(t+\tau)\right)_Y\,dt=-\left(a,W^{0}\dot
f_k^{L_k}(\cdot+\tau)\right)_H=\left(a,\overline{\lambda}_kW^0f_k^{L_k}(\cdot+\tau)\right)_H=\notag\\
\left(\lambda_k a,W^{0}f_k^{L_k}(\cdot+\tau)\right)_H=-\int_0^{2T}
\left(\lambda_k(Oa)(t),f_k^{L_k}(t+\tau)\right)_Y\,dt\label{DV}
\end{eqnarray}
and for $i<L_k$:
\begin{eqnarray}
\int_0^{2T}\left((Oa)(t),\dot
f_k^{i}(t+\tau)\right)_Y\,dt=\left(a,\overline{\lambda}_kW^0f_k^{i}(\cdot+\tau)+W^0f_k^{i+1}(\cdot+\tau)\right)_H=\notag\\
-\lambda_k\int_0^{2T}
\left((Oa)(t),f_k^{i}(t+\tau)\right)_Y\,dt-\int_0^{2T}
\left((Oa)(t),f_k^{i+1}(t+\tau)\right)_Y\,dt\label{DV1}
\end{eqnarray}

In what follows we  assume that elements with index
$i=L_k+1$ or $i=0$ are zero. Combining (\ref{OD}) and (\ref{DV}),
(\ref{DV1}), we see that the pair $\lambda_k,\, f_k$ satisfies on
$0<\tau<T$, $i=1,\ldots, L_k$:
\begin{equation}
\label{F_eqn}
\int_0^{2T}\left(\dot{(Oa)}(t)-\lambda_k(Oa)(t),f_k^i(t+\tau)\right)_Y\,dt=\int_0^{2T}\left((Oa)(t),f_k^{i+1}(t+\tau)\right)_Y\,dt.
\end{equation}

 Now we prove the converse:  solving the generalized
eigenvalue problem
\begin{equation}
\label{F_eqn1}
\int_0^{2T}\left(\dot{(Oa)}(t)-\lambda(Oa)(t),f(t+\tau)\right)_Y\,dt=0
\end{equation}
yields $\{\lambda_k\}_{k=1}^\infty$ eigenvalues of $A$ and
controls $\{f_k^i\}$, $i=1,\ldots,L_k,$  $k \in \mathbb{N}.$

Let the functions $\{f_1,\ldots,f_L\}$ satisfying (\ref{F_eqn})
constitute the chain for (\ref{F_eqn1}) for some $\lambda$. Then
as it follows from the proof that for $\tau\in (0,T)$:
\begin{equation*}
\left(a,W^{0}\dot
f_i(t+\tau)\right)_H+\lambda\left(a,W^0f_i(t+\tau)\right)_H=-\left(a,W^0f_{i+1}(t+\tau)\right)_H,
\end{equation*}
which is equivalent to
\begin{equation}
\label{rev_2}
-\left(a,A^*v^{f_i(t+\tau)}(0)\right)_H+\lambda\left(a,v^{f_i(t+\tau)}(0)\right)_H=-\left(a,v^{f_{i+1}(t+\tau)}(0)\right)_H,\,\,
\tau\in (0,T).
\end{equation}
First we consider case $i=L$. Rewriting the last equality (we use
the notation $f=f_L$) as
\begin{equation}
\label{rev_3} \left(a,A^*v^{f(t+\tau)}(0)-\overline\lambda
v^{f(t+\tau)}(0)\right)_H=0,\quad \tau\in (0,T).
\end{equation}
We assume that
$v^{f(t+\tau)}(T-\tau)=\sum_{{k \in \mathbb{N}}\atop{i=1,\ldots,
L_k}} c_k^i\psi_k^i$. Then developing $v^f$ in the Fourier series as
we did in (\ref{v_f}), we arrive at:
\begin{equation}
\label{vf}
v^{f(t+\tau)}(0)=\sum_{{k\in \mathbb{N}}\atop{i=1,\ldots, L_k}}
c_k^i\sum_{j=1}^{L_k}\frac{(T-\tau)^{j-i}}{(j-i)!}e^{\overline{\lambda}_k(T-\tau)}\psi_k^j
\end{equation}
Applying operator $A^*$ and using the property
$A^*\psi_k^j=\overline{\lambda}_k\psi_k^j+\psi_k^{j+1},$ we
obtain:
\begin{equation}
\label{Af}
A^*v^{f(t+\tau)}(0)=\sum_{{k\in \mathbb{N}}\atop{i=1,\ldots,
L_k}}
c_k^i\sum_{j=1}^{L_k}\frac{(T-\tau)^{j-i}}{(j-i)!}e^{\overline{\lambda}_k(T-\tau)}\left(\overline{\lambda}_k\psi_k^j+\psi_k^{j+1}\right)
\end{equation}
Introducing the notation
\begin{equation}
\label{g} g(\tau):=A^*v^{f(t+\tau)}(0)-\overline\lambda
v^{f(t+\tau)}(0)=\sum_{{k \in \mathbb{N}}\atop{i=1,\ldots,
L_k}}g_k^i(\tau)\psi_k^i,
\end{equation}
relation (\ref{rev_3}) yields:
\begin{equation}
\label{eq_zero}
0=\left(a,g\right)_H=\sum_{{k\in \mathbb{N}}\atop{i=1,\ldots,
L_k}}a_k^i g_k^i(\tau),\quad \tau\in (0,T).
\end{equation}
The functions $g_k^i(\tau)$ are combination of products of
$e^{\overline{\lambda}_k(T-\tau)}$ and polynomials
$\frac{(T-\tau)^\alpha}{\alpha !}$. Then we can rewrite
(\ref{eq_zero})  as
\begin{equation}
\label{eq_zero1} 0=\sum_{{k\in \mathbb{N}}\atop{i=1,\ldots,
L_k}}b_k^i
\frac{(T-\tau)^{i-1}}{(i-1)!}e^{\overline{\lambda}_k(T-\tau)},\quad
\tau\in (0,T).
\end{equation}

If $Y=\mathbb{R}$,
the controllability of the dynamical system
(\ref{E_AD}) imply  \cite{AI} the minimality of the family
$\bigcup_{k=1}^\infty\{e^{\overline{\lambda}_kt},te^{\overline{\lambda}_kt}\ldots,t^{L_k-1}e^{\overline{\lambda}_kt}\}$
in $L_2(0,T)$ in $L_2(0,T)$, so we have $b_k^i=0$ for all $k,i.$
On the other hand, as follows from (\ref{vf}), (\ref{Af}):
\begin{equation*}
b_k^{L_k}=c_k^1\overline{\lambda}_ka_k^1-\overline{\lambda}c_k^1a_k^1=0.
\end{equation*}
Then since $a$ is generic, either $\lambda=\lambda_k$ or
$c_k^1=0$.

Let $\lambda\not=\lambda_k$, so $c_k^1=0$. Then for $b_k^{L_k-1}$
we have:
\begin{equation*}
b_k^{L_k-1}=c_k^2\overline{\lambda}_ka_k^2-\overline{\lambda}c_k^2a_k^2=0,
\end{equation*}
from which the equality $c_k^2=0$ follows. Repeating this
procedure for $b_k^{L_k-i}$, $i\geqslant 2$, we obtain:
\begin{equation}
\label{cond1}
\text{If $\lambda\not=\lambda_k$, then $c_k^i=0,$
$i=1,\ldots,L_k$}.
\end{equation}

Consider the second option: let $\lambda=\lambda_k$. Then from
(\ref{vf}), (\ref{Af}):
\begin{equation*}
b_k^{L_k-1}=c_k^1=0,\quad b_k^{L_k-2}=c_k^2a_k^3=0,\ldots,
b_k^{1}=c_k^{L_k-1}a_k^{L_k}=0.
\end{equation*}
So we arrive at
\begin{equation}
\label{cond2} \text{If $\lambda=\lambda_k$, then $c_k^i=0$,
$i=1,\ldots, L_k-1$, and $c_k^{L_k}$ could be arbitrary.}
\end{equation}
Finally (\ref{cond1}), (\ref{cond2}) imply that
$\lambda=\lambda_{k'}$ and $f=c_{k'}f_{k'}^{L_{k'}}$,
$c_{k'}\not=0$, for some $k'.$

 Thus on the first step we already obtained that
$\lambda=\lambda_{k'}$ for some $k'$ and
$f_L=c_{k'}f_{k'}^{L_{k'}}$. The second vector $f$ in the Jordan
chain satisfies
\begin{equation*}
\int_0^{2T}\left(\dot{(Oa)}(t)-\lambda_{k'}(Oa)(t),f(t+\tau)\right)_Y\,dt=\int_0^{2T}\left((Oa)(t),c_{k'}f_{k'}^{L_{k'}}(t+\tau)\right)_Y\,dt.
\end{equation*}
We rewrite (\ref{rev_2}) in our case:
\begin{equation}
\label{rev_4}
-\left(a,A^*v^{f(t+\tau)}(0)\right)_H+\lambda_{k'}\left(a,v^{f(t+\tau)}(0)\right)_H=-\left(a,c_{k'}v^{f_{k'}^{L_{k'}}(t+\tau)}(0)\right)_H,\,\,
\tau\in (0,T).
\end{equation}
In this case $g$ introduced in (\ref{g}) has a form
\begin{equation*}
g(\tau)=\sum_{{k \in \mathbb{N}}\atop{i=1,\ldots, L_k}}
c_k^i\sum_{j=1}^{L_k}\frac{(T-\tau)^{j-i}}{(j-i)!}e^{\overline{\lambda}_k(T-\tau)}\left(\left(\overline{\lambda}_k-\overline\lambda_{k'}\right)\psi_k^j+\psi_k^{j+1}\right)
\end{equation*}
and rewrite (\ref{rev_4}) as
\begin{equation}
\label{eq_zero2}
\left(a,g\right)_H=\left(a,v^{f_{k'}^{L_{k'}}(\cdot+\tau)}\right)_H=c_{k'}a_{k'}^{L_{k'}}e^{\overline{\lambda}_{k'}(T-\tau)}
\end{equation}
Using the same notations as for (\ref{eq_zero}), (\ref{eq_zero1}),
we write down the equalities for coefficients $b_k^i$ for
(\ref{eq_zero2}) to get:
\begin{equation*}
b_{k'}^1=c_{k'}a_{k'}^{L_{k'}},\quad b_k^i=0, \quad k\not= k',
\end{equation*}
In the case $k\not=k'$ we repeat the arguments used above and find
that
$$
c_k^i=0,\quad i=1,\ldots,L_k.
$$
When $k=k'$, we have:
\begin{eqnarray*}
b_{k'}^{L_{k'}}=0,\quad
b_{k'}^{L_{k'}-1}=c_{k'}^1a_{k'}^{L_{k'}}=0,\quad
b_{k'}^{L_{k'}-2}=c_{k'}^2a_{k'}^{L_{k'}}=0,\\
b_{k'}^2=c_{k'}^{L_{k'}-2}a_{k'}^{L_{k'}}=0,\quad
b_{k'}^1=c_{k'}^{L_{k'}-1}a_{k'}^{L_{k'}}=c_{k'}a_{k'}^{L_{k'}}.
\end{eqnarray*}
So we find:
\begin{equation*}
c_{k'}^i=0, \quad i< L_{k'}-1,\quad
c_{k'}^{L_{k'}-1}=c_{k'},\quad\text{$c_{k'}^{L_{k'}}$ is
arbitrary}
\end{equation*}
So finally we arrive at for some $c_{L-1}$:
\begin{equation*}
f=f_{L-1}=c_{k'}f_{k'}^{L_{k'}-1}+c_{L-1}f_{k'}^{L_{k'}}
\end{equation*}
Arguing in the same fashion, we obtain that
\begin{equation*}
f_i=c_{k'}f_{k'}^{L_{k'}-i}+c_if_{k'}^{L_{k'}},\quad 1\leqslant i<
L_{k'}-1.
\end{equation*}
So we have shown that the elements of the Jordan chain for
(\ref{M_eqn}) which correspond to eigenvalue $\lambda_{k'}$  are
the sum of corresponding controls and eigenvector (i.e. the
control that generate the eigenvector of $A^*$).
\end{proof}
\begin{remark}
\label{remark}
 The solution to (\ref{M_eqn}) yields
$\{\lambda_k\}_{k=1}^\infty$ eigenvalues of $A$ and
(non-normalized) root vectors $\{\widehat f_k^i\}$, $\widehat
f_k^i= c_{k}f_k^i+c_k^i f_k^{L_k}$ $k\in \mathbb{N}$,
$i=1,\ldots,L_k,$ $c_k^{L_k}=0$.
\end{remark}

For the dynamical system (\ref{E_D}), under the conditions on $A,$
$Y$, formulated in Theorem 1, there is the possibility to extend
the observation $y(t)=\left(Ou^a\right)(t)$  defined for $t\in
(0,2T)$ to $t\in \mathbb{R}_+$. To this aim we show that for
observation having a form
\begin{equation}
\label{Observ} \mathbb{O}a=\sum_{k\in \mathbb{N}}
e^{\lambda_kt}\sum_{j=1}^{L_k}\frac{b_k^jt^{L_k-j}}{(L_k-j)!}
\end{equation}
we can recover the coefficients $b_k^j.$

Take $i\in \{1,\ldots,L_k\}$ and look for the solution to
(\ref{E_D}) with $a=\phi_k^i$ in the form
$u=\sum_{l=1}^{L_k}c_l(t)\phi_k^l$, we arrive at the system (here
$c_{L_k+1}=0$):
\begin{align*}
\frac{d}{dt}c_l(t)-\lambda_kc_l(t)=c_{l+1}(t),\quad
l=1,\ldots,L_k,\\
c_l(0)=\delta_{li}.
\end{align*}
whose solution is
\begin{eqnarray*}
c_l(t)=\frac{t^{i-l}}{(i-l)!}e^{\lambda_kt},\quad l\leqslant i,\\
c_l(t)=0,\quad l>i.
\end{eqnarray*}
Thus
\begin{equation}
\label{u_phi}
u^{\phi_k^i}=\sum_{l=1}^i\frac{t^{i-l}}{(i-l)!}e^{\lambda_kt}\phi_k^l.
\end{equation}
For the initial state $a=\sum_{k\in
\mathbb{N}}\sum_{i=1}^{L_k}a_k^i\phi_k^i$ we obtain:
\begin{equation*}
u^{a}=\sum_{k\in
\mathbb{N}}e^{\lambda_kt}\sum_{j=1}^{L_k}\frac{t^{L_k-j}}{(L_k-j)!}\sum_{l=1}^j
a_k^{L_k-j+l}\phi_k^l.
\end{equation*}
So for observation $(\mathbb{O}a)(t)=\left(Ou^a\right)(t)$ we get
the representation (\ref{Observ}) with coefficients $b_k^j$
defined by
\begin{equation}
\label{bj} b_k^j:=\sum_{l=1}^j a_k^{L_k-j+l}O\phi_k^l, \quad k\in
\mathbb{N},\,\, j=1,\ldots, L_k.
\end{equation}
Making use of Theorem \ref{teor} (see also Remark \ref{remark}),
we have:
\begin{equation}
\label{Control} W^0\widehat
f^i_k=c_k\psi_k^i+c_k^i\psi_k^{L_k},\quad k\in \mathbb{N},\,\,
i=1,\ldots,L_k,\,\,c_k^{L_k}=0.
\end{equation}
Counting (\ref{Adj}) we write:
\begin{equation*}
\left(W^0\widehat f_k^i,a\right)_H=-\int_0^T Ou^a\widehat
f_k^i\,dt.
\end{equation*}
We plug $a=\phi_k^i$ in the last equality and use (\ref{Control})
to get
\begin{equation}
\label{c_k}
c_k=\left(c_k\psi_k^i+c_k^i\psi_k^{L_k},\phi_k^i\right)_H=-\int_0^T
Ou^{\phi_k^i}\widehat f_k^i\,dt.
\end{equation}
We evaluate the right hand side of (\ref{c_k}) for all $i$. For
$i=1$ we get (see (\ref{u_phi})):
\begin{equation*}
c_k=-O{\phi_k^1}\int_0^T e^{\lambda_kt}\widehat f_k^1\,dt.
\end{equation*}
Or equivalently:
\begin{equation}
\label{c_k1} \frac{c_k}{O{\phi_k^1}}=-\int_0^T
e^{\lambda_kt}\widehat f_k^1\,dt.
\end{equation}
Evaluating (\ref{c_k}) for $i=2$, counting (\ref{u_phi}), we
obtain:
\begin{equation*}
c_k=-O{\phi_k^2}\int_0^T e^{\lambda_kt}\widehat
f_k^2\,dt-O{\phi_k^1}\int_0^T te^{\lambda_kt}\widehat f_k^2\,dt.
\end{equation*}
Divide this equality by $c_k$ and plug (\ref{c_k1}) to find:
\begin{equation}
\label{c_k2} \frac{c_k}{O{\phi_k^2}}=-\frac{\int_0^T
e^{\lambda_kt}\widehat f_k^1\,dt \int_0^T e^{\lambda_kt}\widehat
f_k^2\,dt}{\int_0^T e^{\lambda_kt}\widehat f_k^1\,dt-\int_0^T
te^{\lambda_kt}\widehat f_k^2\,dt}
\end{equation}
Suppose we already found $\frac{c_k}{O{\phi_k^l}}$ for
$l=1,\ldots,i-1$. To find this quantity for $l=i$, we evaluate
(\ref{c_k}), plugging expression for $u^{\phi_k^i}$ (\ref{u_phi}):
\begin{equation*}
c_k=-\sum_{l=1}^i O{\phi_k^l}\int_0^T \frac{t^{i-l}}{(i-l)!}
e^{\lambda_kt}\widehat f_k^i\,dt.
\end{equation*}
We divide last equality by $c_k$ to find:
\begin{equation}
\label{c_ki} \frac{c_k}{O{\phi_k^i}}=-\frac{\int_0^T
e^{\lambda_kt}\widehat f_k^i\,dt}{1+\sum_{l=1}^{i-1} \int_0^T
\frac{t^{i-l}}{(i-l)!} e^{\lambda_kt}\widehat
f_k^i\,dt\left(\frac{c_k}{O{\phi_k^l}}\right)^{-1}}
\end{equation}
Observe that in the right hand side of (\ref{c_ki}) in view of
(\ref{c_k2}), we know all terms.

To evaluate $a_k^i$ we use, see (\ref{Control}):
\begin{equation}
\label{a_ki} a_k^i=\left(a,\psi_k^i\right)_H=\left(a,W^0\widehat
f_k^i-c_k^i\psi_k^{L_k}\right)_H\frac{1}{c_k}=-\int_0^T
Ou^a\widehat f_k^i\,dt\frac{1}{c_k}-a_k^{L_k}\frac{c_k^i}{c_k}
\end{equation}
We multiply (\ref{Control}) by  $\phi_k^{L_k}$ and get for
$i<L_k$:
\begin{eqnarray*}
c_k^i=\left(W^0f_k^i,\phi_k^{L_k}\right)_{H}=-\int_0^T
f_k^i(t)\left(Ou^{\phi_k^{L_k}}\right)(t)\,dt\\
=-\sum_{l=1}^{L_k}O\phi_k^l\int_0^T\frac{t^{L_k-l}}{(L_k-l)!}e^{\lambda_kt}f_k^i(t)\,dt.
\end{eqnarray*}
Dividing the last equality by $c_k$ we get
\begin{equation}
\label{c_kratio}
\frac{c_k^i}{c_k}=-\sum_{l=1}^{L_k}\left(\frac{c_k}{O\phi_k^l}\right)^{-1}\int_0^T\frac{t^{L_k-l}}{(L_k-l)!}e^{\lambda_kt}f_k^i(t)\,dt,\quad
i<L_k.
\end{equation}
Notice that in view of (\ref{c_ki}), we know all terms in the
right hand side in (\ref{c_kratio}). Now we multiply (\ref{a_ki})
by $c_k$:
\begin{equation*}
a_k^ic_k=-\int_0^T Ou^a\widehat
f_k^i\,dt-a_k^{L_k}c_k\frac{c_k^i}{c_k}.
\end{equation*}
Since $c_k^{L_k}=0,$ we have for $i=L_k$:
\begin{equation*}
a_k^{L_k}c_k=-\int_0^T\widehat f_k^{L_k}(t)
\left(Ou^a\right)(t)\,dt,
\end{equation*}
and finally
\begin{equation}
\label{akck} a_k^ic_k=-\int_0^T \widehat
f_k^i(t)\left(Ou^a\right)(t)\,dt+\int_0^T\widehat f_k^{L_k}(t)
\left(Ou^a\right)(t)\,dt\frac{c_k^i}{c_k}.
\end{equation}
In view of (\ref{c_kratio}), we know all terms in the right hand
side of (\ref{akck}).

Now we rewrite formula for $b_k^j$ (\ref{bj}):
\begin{equation}
\label{bjnew} b_k^j:=\sum_{l=1}^j
\left\{a_k^{L_k-j+l}c_k\right\}\left(\frac{O\phi_k^l}{c_k}\right)\quad
k\in \mathbb{N},\,\, j=1,\ldots, L_k.
\end{equation}
and observe that the first term in each summand is given by
(\ref{akck}), while the second term by (\ref{c_ki}). So we know
right hand side in (\ref{bjnew}).

After we recovered all $b_k^j$ by (\ref{bjnew}), we can extend the
observation $(\mathbb{O}a)(t)$ by formula (\ref{Observ}) for
$t>2T$.

\section{Application to inverse problems.}

Here provide two application of the theory developed above to
inverse problems. Other applications of the BC approach to the
spectral estimation problem
 can be found in \cite{AB,ABN,AGM,AM1,AM2,AMR,MM}.

\subsection{Reconstructing the potential for the 1D Schr\"odinger
equation from boundary measurements}

Let the real potential $q\in L^1(0,1)$ and $a\in H^1_0(0,1)$ be
fixed, we consider the boundary value problem:
\begin{equation}
    \label{eq_1}
    \left\{
    \begin{array}{ll}
         iu_t(x,t)-u_{xx}(x,t)+q(x)u(x,t)=0 &\qquad t>0,\quad 0<x<1\\
         u(0,t)=u(1,t)=0&\qquad t>0,\\
         u(x,0)=a(x) & \qquad 0<x<1.
    \end{array}
    \right.
\end{equation}
Assuming that the initial datum $a$ is generic (but unknown), the
inverse problem we are interested in is to determine the potential
$q$ from the trace of the derivative of the solution $u$ to
(\ref{eq_1}) on the boundary:
$$\{r_0(t),r_1(t)\}:=\{u_x(0,t),u_x(1,t)\}, \quad t\in (0,2T),$$
It is well known that the selfadjoint operator $A$ defined  on
$L^2(0,1)$  by
\begin{equation}\label{eqA}
A\phi = -\phi'' +q\phi,\quad D(A):=H^2(0,1)\cap H^1_0(0,1).
\end{equation}
admits a family of eigenfunctions $\{\phi_k\}_{k=1}^\infty$
forming a orthonormal basis in $L^2(0,1)$, and associated sequence
of eigenvalues $\lambda_k\rightarrow +\infty$. Using the Fourier
method, we can represent the solution of (\ref{eq_1}) in the form
\begin{equation}
\label{Repr} u(x,t)=\sum_{k=1}^\infty a_ke^{i\lambda_k
t}\phi_k(x),\quad a_k=\left(a,\phi_k\right)_{L^2(0,1)}
\end{equation}
The inverse data admits the representation
\begin{equation}
\left\{r_0(t),r_1(t)\right\}=\left\{\sum_{k=1}^\infty
a_ke^{i\lambda_k t}\phi_k'(0),\sum_{k=1}^\infty a_ke^{i\lambda_k
t}\phi_k'(1)\right\}.\label{RT}
\end{equation}
One can prove that $r_0, r_1\in L^2(0,T)$. Using the method from
the first section, we recover the eigenvalues $\lambda_k$ of $A$
and the products $\phi_k'(0)a_k$ and $\phi_k'(1)a_k$. So (as $a$
is generic) we recovered the spectral data consisting of
\begin{equation}
\label{Data}
D:=\left\{\lambda_k,\frac{\phi_k'(1)}{\phi_k'(0)}\right\}_{k=1}^\infty.
\end{equation}
Now from $D$ we construct the spectral function associated to $A$.

Given $\lambda\in {\mathbb C}$, denote by $y(\cdot,\lambda)$ the
solution to
\begin{equation*}
\left\{
\begin{array}l
- y''(x,\lambda) + q(x)y(x,\lambda)=\lambda y(x,\lambda),\qquad
0<x<1, \\
y(0,\lambda)=0,\quad y'(0,\lambda)=1.
\end{array}
\right.
\end{equation*}
Then the eigenvalues of the Dirichlet problem  of $A$ are exactly
the zeros of the function $y(1,\lambda)$, while a family of
normalized corresponding eigenfunctions is given by
$\phi_k(x)=\dfrac{y(x,\lambda_k)}{\|y(\cdot,\lambda_k)\|}$. Thus
we can rewrite the second components in $D$ in the following way:
\begin{equation}
\label{A_k}
\frac{\phi_k'(1)}{\phi_k'(0)}=\frac{y'(1,\lambda_k)}{y'(0,\lambda_k)}=y'(1,\lambda_k)=:A_k.
\end{equation}
Let us denote by dot the derivative with respect to $\lambda$ and
$\lambda_n$ be an eigenvalue of $A$. We borrowed the following
fact from \cite[p. 30]{PT}:
\begin{align*}
\|y(\cdot,\lambda_k)\|^2_{L^2}=y'(1,\lambda_k)\dot
y(1,\lambda_k),\\
y(1,\lambda)=\prod_{n\geqslant
1}\frac{\lambda_n-\lambda}{n^2\pi^2}\\
\dot y(1,\lambda_k)=-\frac{1}{k^2\pi^2}\prod_{n\geqslant 1,
n\not=k}\frac{\lambda_n-\lambda_k}{n^2\pi^2}=:B_k.
\end{align*}
Notice that the set of pairs
$\{\lambda_k,\|y(\cdot,\lambda_k)\|^2_{L^2}\}_{k=1}^\infty=:\widetilde
D$ is a ``classical'' spectral data. Using the above relations, we
come to $\widetilde
D=\left\{\lambda_k,A_kB_k\right\}_{k=1}^\infty$. Let
$\alpha_k^2:=\|y(\cdot,\lambda_k)\|^2_{L^2}=A_kB_k$, we introduce
the spectral function associated with $A$:
\begin{equation*}
\label{sp_func} \rho(\lambda)=\left\{\begin{array}l
-\sum\limits_{\lambda\leqslant\lambda_k\leqslant
0}\frac{1}{\alpha_k^2}\quad \lambda\leqslant 0, \\
\sum\limits_{0<\lambda_k\leqslant\lambda}\frac{1}{\alpha_k^2}\quad
\lambda
> 0,
\end{array}\right.
\end{equation*}
which is a monotone increasing function having jumps at the points
of the Dirichlet spectra. The regularized spectral function is
introduced by
\begin{equation*}
\label{reg_sp_func}
\sigma(\lambda)=\left\{\begin{array}l \rho(\lambda)-\rho_0(\lambda)\quad \lambda\geqslant 0, \\
\rho(\lambda)\quad \lambda < 0,
\end{array}\right. \quad
\rho_0(\lambda)=\sum\limits_{0<\lambda_k^0\leqslant\lambda}\frac{1}{(\alpha_k^0)^2}\quad
\lambda
> 0,
\end{equation*}
where $\rho_0$ is the spectral function associated with the
operator $A$ with $q\equiv 0$. The potential can thus be recovered
from $\sigma(\lambda)$ by Gelfand-Levitan, Krein or the BC method
(see \cite{AM,BM_1}). Once the potential has been found, we can
recover the eigenfunctions $\phi_k$, the traces $\phi_k'(0)$ and
Fourier coefficients $a_k$, $k=1,\ldots \infty$. Thus, the initial
state can be recovered via its Fourier series.

\subsection{Extension of the inverse data}

We fix $p_{ij}\in C^1([0,1];\mathbb{C})$, $d_1, d_2 \in
L_2(0,1;\mathbb{C})$ and consider on interval $(0,1)$ the initial
boundary value problem
\begin{equation}
\label{eq} \left\{\begin{array}l
\frac{\p}{\p t}\begin{pmatrix} u\\
v\end{pmatrix}-\frac{\p}{\p x}\begin{pmatrix} 0&1\\
1&0\end{pmatrix}\begin{pmatrix} u\\
v\end{pmatrix}-\begin{pmatrix} p_{11}& p_{12}\\
p_{21}&p_{22}\end{pmatrix}\begin{pmatrix} u\\
v\end{pmatrix}=0, \quad t>0, \\
u(0,t)=u(1,t)=0,\quad t>0,\\
\begin{pmatrix} u(x,0)\\
v(x,0)\end{pmatrix}=\begin{pmatrix} d_1(x)\\
d_2(x)\end{pmatrix},\quad 0\leqslant x\leqslant 1
\end{array}
\right.
\end{equation}
We fix some $T>0$ and define $R(t):=\{v(0,t),v(1,t)\}$,
$0\leqslant t\leqslant T.$ Here we focus on the problem of the
continuation of the inverse data: we assume that $R(t)$ is known
on the interval $(0,T)$, $T>2$, and recover it on the whole real
axis. The problem of the recovering unknown coefficients $p_{ij}$
and initial state $c_{1,2}$ has been considered in \cite{TY1,TY5},
where the authors established the uniqueness result, having the
response $R(t)$ on the interval $(-T,T)$ for large enough $T.$

We introduce the notations $B=\begin{pmatrix} 0&1\\
1&0\end{pmatrix}$, $P=\begin{pmatrix} p_{11}& p_{12}\\
p_{21}&p_{22}\end{pmatrix}$, $D=\begin{pmatrix} d_1\\
d_2\end{pmatrix}$  and the operators $A,$ $A^*$ acting by the rule
\begin{eqnarray*}
A=B\frac{d}{dx}+P,\quad \text{on $(0,1)$},\\
A^*\psi=-B\frac{d}{dx}+P^T,\quad \text{on $(0,1)$},
\end{eqnarray*}
with the domains
$$
D(A)=D(A^*)=\left\{\varphi=\begin{pmatrix} \varphi_1\\
\varphi_2\end{pmatrix}\in H^1(0,1;\mathbb{C}^2)\,|\,
\varphi_1(0)=\varphi_1(1)=0\right\}
$$
The spectrum of the operator $A$ has the following structure (see
\cite{TY1,TY5}): $\sigma(A)=\Sigma_1\cup\Sigma_2,$ where
$\Sigma_1\cap\Sigma_2=\emptyset$ and there exists $N_1 \in \N$
such that
\begin{itemize}
\item[1)]$\Sigma_1$ consists of $2N_1-1$ eigenvalues including
algebraical multiplicities

\item[2)] $\Sigma_2$ consists of infinite number of eigenvalues of
multiplicity one

\item[3)] Root vectors of $A$ form  a Riesz basis in $
L_2(0,1;\mathbb{C}^2)$.
\end{itemize}

Let  $m$  denote the algebraical multiplicity of eigenvalue
$\lambda,$ and we introduce the notations:
\begin{eqnarray*}
\Sigma_1=\left\{\lambda^i\in\sigma(A), \, m_i\geqslant 2,\,
1\leqslant i\leqslant N\right\},\\
\Sigma_2=\left\{\lambda_n\in\sigma(A), \, \lambda_n \text{ is
simple },\, n\in \mathbb{Z}\right\}.
\end{eqnarray*}
Let $e_1:=\left(0\atop 1\right)$. The root vectors are introduced
in the following way:
\begin{eqnarray*}
\left(A-\lambda^i\right)\phi^i_1 =0,\quad
\left(A-\lambda^i\right)\phi^i_j =\phi^i_{j-1},\quad 2\leqslant
j\leqslant m_i,\\
\phi^i_j(0)=e_1, \,\, \phi^i_j\in D(A),\,\, 1\leqslant j\leqslant
m_i.
\end{eqnarray*}
For the adjoint operator the following equalities are valid:
\begin{eqnarray*}
\left(A^*-\overline\lambda^i\right)\psi^{i}_{m_i} =0,\quad
\left(A^*-\overline\lambda^i\right)\psi^{i}_j
=\psi^{i}_{j+1},\quad 1\leqslant
j\leqslant m_i-1,\\
\psi^{i}_j(0)=e_1, \,\, \psi^{i}_j\in D(A^*),\,\, 1\leqslant
j\leqslant m_i.
\end{eqnarray*}
For the simple eigenvalues we have:
\begin{eqnarray*}
\left(A-\lambda_n\right)\phi_n =0,\quad
\left(A^*-\overline\lambda_n\right)\psi_{n} =0,\\
\phi_n(0)=\psi_{n}(0)=e_1, \,\, \phi_{n}\in D(A),\,\,\psi_{n}\in
D(A^*).
\end{eqnarray*}
Moreover, the following biorthogonality conditions hold:
\begin{align*}
\left(\phi^i_j,\psi_{n}\right)=0,\quad
\left(\phi_n,\psi^{i}_j\right)=0,\quad
\left(\phi_k,\psi_{n}\right)=0,\\
\left(\phi^i_j,\psi^{k}_l\right)=0,\quad \text{if $i\not= k$ or
$j\not= l$},\\
\rho^i_j=\left(\phi^i_j,\psi^{i}_j\right),\quad i=1,\ldots,N,\quad
j=1,\ldots, m_i,\\
\rho_n=\left(\phi_n,\psi_{n}\right),\quad n\in \mathbb{Z},
\end{align*}
We represent the initial state as the series:
\begin{equation}
\label{D_repr} D=\sum_{i=1}^N
\sum_{j=1}^{m_i}d^i_j\phi^i_j(x)+\sum_{n\in
\mathbb{Z}}d_n\phi_n(x).
\end{equation}
and look for the solution to (\ref{eq}) in the form
\begin{equation*}
\begin{pmatrix}u\\ v\end{pmatrix}(x,t)=\sum_{i=1}^N
\sum_{j=1}^{m_i}c_j^i(t)\phi_j^i(x)+\sum_{n\in
\mathbb{Z}}c_n(t)\phi_n(x).
\end{equation*}
Using the method of moments we can derive the system of ODe's for
$c^i_j,$ $i\in\{1,\ldots,N\}$, $j\in \{1,\ldots,m_i\}$; $c_n$,
$n\in \mathbb{Z}$ solving which we obtain
\begin{eqnarray*}
c^i_j(t)=e^{\lambda^i t}\left[d^i_j+d^i_{j+1}t+d^i_{j+2}\frac{t^2}{2}+\ldots+d^i_{m_i}\frac{t^{m_i-j}}{(m_i-j)!}\right],\\
c_n(t)=d_ne^{\lambda_n t}.
\end{eqnarray*}
Notice that the response $\{v(0,t),v(1,t)\}$ has a form depicted
in (\ref{signal}):
\begin{eqnarray}
v(0,t)=\sum_{i=1}^N e^{\lambda^i t} a^0_i(t)+\sum_{n\in \mathbb{Z}}e^{\lambda_n t}d_n\left(\phi_n(0)\right)_2,\label{r1}\\
v(1,t)=\sum_{i=1}^N e^{\lambda^i t} a^1_i(t)+\sum_{n\in
\mathbb{Z}}e^{\lambda_n t}d_n\left(\phi_n(1)\right)_2,\label{r2}
\end{eqnarray}
where the coefficients of
$a^0_i(t)=\sum_{k=0}^{m_i-1}\alpha_k^it^k$ are given by
\begin{eqnarray*}
\alpha^i_0=\sum_{l=1}^{m_i}d^i_l\left(\phi^i_{l}(0)\right)_2,\quad \alpha^i_1=\sum_{l=2}^{m_i}d^i_l\left(\phi^i_{l-1}(0)\right)_2,\quad \alpha^i_2=\frac{1}{2}\sum_{l=3}^{m_i}d^i_l\left(\phi^i_{l-2}(0)\right)_2,\\
\ldots,
\alpha^i_k=\frac{1}{k!}\sum_{l=k+1}^{m_i}d^i_l\left(\phi^i_{l-k}(0)\right)_2,\ldots\quad
\alpha^i_{m_i-1}=\frac{1}{(m_i-1)!}d^i_{m_i}\left(\phi^i_{1}(0)\right)_2.
\end{eqnarray*}
The coefficients $a^1_i(t)$, $i=1,\ldots,N$ are defined by the
similar formulas.

We assume that the initial state $D$ is generic. Introducing the
notation $U:= \left(u\atop v\right)$ we consider the dynamical
system with the boundary control $f\in L_2(\mathbb{R}_+)$:
\begin{equation*}
\left\{
\begin{array}l
U_t-AU=0, \quad 0\leqslant x\leqslant 1,\,\, t>0,\\
u(0,t)=f(t), u(1,t)=0,\quad  t>0, \\
U(x,0)=0.
\end{array}
\right.
\end{equation*}
It is not difficult to show that this system is exactly
controllable in time $T \geq 2$. This implies (see \cite{AI}) that
the family
$\bigcup_{i=1}^N\{e^{\lambda^it},\ldots,t^{m_i-1}e^{\lambda^it}\}\cup\{e^{\lambda_nt}\}_{n\in
\mathbb{Z}}$ forms a Riesz basis in a closure of its linear span
in $L_2((0,T);\mathbb{C})$. So we can apply the method from the
second sections to recover $\lambda^i$, $m_i$, coefficients of
polynomials $a^{0,1}_i(t)$ $i=1,\ldots,N,$ $\lambda_n$, $n\in
\mathbb{Z}$. The latter allows one to extend the inverse data
$R(t)$ to all values of $t\in \mathbb{R}$ by formulas (\ref{r1}),
(\ref{r2}). This is important to the solution of the
identification problem, see \cite{TY5}.

\noindent{\bf Acknowledgments} Sergei Avdonin  was supported by
the NSF grant DMS 1411564; Alexander Mikhaylov was supported by
NSh-1771.2014.1, RFBR 14-01-00306; Victor Mikhaylov was supported
by RFBR 14-01-00535, RFBR 14-01-31388 and NIR SPbGU
11.38.263.2014.

\end{document}